\documentclass[a4paper,12pt,fleqn]{article}

\usepackage{amssymb}
\usepackage{amsmath}
\usepackage{graphicx}
\usepackage{color}
\usepackage{eurosym}
\usepackage{txfonts}
\usepackage{courier}
\usepackage{pifont}

\usepackage[colorlinks=false,urlcolor=blue,ps2pdf=true,
pdfpagemode=UseNone, pdfauthor={Robin Whitty},
pdfsubject={Triangle bisection}, pdfcreator={Robin Whitty},
pdfproducer={Robin Whitty}]{hyperref}

\newcommand{\vc}[1]{{ \normalfont {\bf #1}}}
\newcommand{\vcc}[1]{\overrightarrow{#1}}

\newcommand{\vo}[1]{{ \normalfont {\bf #1}}^{\bot}}

\newcommand{\matctwoone}[2]{
\left(\begin{array}{c} #1 \\ #2
\end{array}\right) }

\author{ Robin
Whitty\thanks{www.theoremoftheday.org}}
\date{\today}
\title{Bisecting a triangle in a given direction}

\begin{document}

\maketitle

\hspace{.25in}
\parbox[c]{6in}{\small {\bf Abstract} Given a triangle, what is the equation of the line which bisects its area and has a given slope? The set of all lines bisecting the area of a triangle has been elegantly determined as a certain `deltoid' envelope and this gives an indirect method of solution. We find that vector algebra allows the equation to be written down rather directly and neatly.
}


Suppose that a triangle is specified by the coordinates of its vertices. A slope is given, and it is required to write down the equation of the line having this slope and bisecting the area of the triangle. In \cite{Dunn}, Dunn and Pretty determine all area bisectors of a triangle in terms of a deltoid-shaped envelope: the line we require will be the unique tangent to this envelope having the given slope. This is very elegant and has inspired some attractive developments, notably in work of Berele and Catoiu: see \cite{Berele}, for example, where much more is recorded about Dunn and Pretty's deltoid. However, it is not the most direct way to approach our problem, particularly because Dunn and Pretty reduce the area bisection question `affinely' to bisecting a right triangle placed at the origin.

\begin{figure}[htb]
\centerline{\includegraphics[scale=0.25]{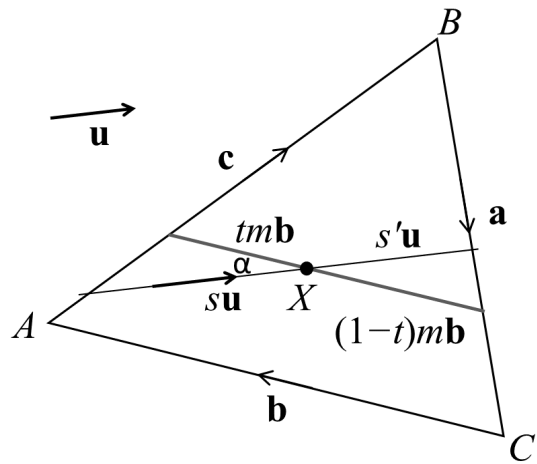}}

\caption{Tilting (anticlockwise) from a bisecting line parallel to edge $\vc{b}$ to create a bisector in the direction of $\vc{u}$.}
\label{fig:vectortriangle}
\end{figure}
We take a different approach, starting with the bisecting line which lies parallel to a chosen triangle edge and tilting it to the required angle. This can be specified quite neatly in terms of vector equations which solve to give that point on the bisecting line through which our desired line must pass. The equation of this line is still not written down in an entirely direct manner because we must choose which triangle edge to start from --- effectively, which Dunn and Pretty deltoid edge we are tangent to --- but this choice turns out to be expressed succinctly in terms of solutions to pairs of simultaneous equations.

Figure~\ref{fig:vectortriangle} shows the vectors which represent our approach algebraically. The  bisecting line parallel to a triangle edge divides the altitude from the opposite vertex in the ratio $1/\sqrt{2}:1-1/\sqrt{2}$. In the figure this is depicted for edge (vector) $\vc{b}$ and we will refer to it as the $\vc{b}$-bisect. Writing $m$ for $1/\sqrt{2}$, the  $\vc{b}$-bisect is  $m\vc{b}$ as a free vector.  The vector whose direction is $\vc{u}$ will correspond to a new bisector of $ABC$ if scalars $t$ and $s$ are such that the triangle with sides $tm\vc{b}$ and $s\vc{u}$ has the same area as the triangle with sides $(1-t)m\vc{b}$ and $s'\vc{u}$, both triangles having the same included angle $\alpha$:
\begin{eqnarray}
\label{eqn:sdash}
 \frac12|s\vc{u}|\times |tm\vc{b}|\sin\alpha &=& \frac12|s'\vc{u}|\times |(1-t)m\vc{b}|\sin\alpha, \nonumber\\
 \mbox{so that  } s'&=&\frac{t}{1-t}s.
\end{eqnarray}
Our aim is to find $t$ since this will specify point $X$ in figure~\ref{fig:vectortriangle}, as a position vector, and hence determine the equation of our bisecting line.

For the moment we are assuming $\vc{b}$ is the right choice of edge direction from which to tilt to give the required direction $\vc{u}$. In fact, our approach will be to specify the collection of all vectors for which~$\vc{b}$ is indeed the right choice. Since for our problem the orientation of these vectors is immaterial we will specify that they follow the orientation of $\vc{b}$ (so, right to left in figure~\ref{fig:vectortriangle}, even though the depicted vector $\vc{u}$ is orientated left to right).

 Now we cannot tilt away from the $\vc{b}$-bisect clockwise  beyond the median line from vertex $C$ to the midpoint of $\vc{c}$, which trivially bisects the triangle area. As a vector this has direction $\vc{b}-\vc{a}$. Nor can we tilt
anticlockwise beyond the median line from the midpoint of $\vc{a}$ to $A$; as a vector this has direction $\vc{b}-\vc{c}$. Then the collection of all free vectors between one median direction or the other, following the orientation of $\vc{b}$,  is given as
$$\vc{u}_b=\vc{b}-\vc{a}+w(\vc{a}-\vc{c}), w\in[0,1].$$
We remark that $w=1/2$ gives the vector $\vc{b}-\frac12(\vc{c}+\vc{a})=\frac32\vc{b}$, precisely the direction of the $\vc{b}$-bisect.

Similarly the collection of vector directions which are valid bisecting tilts from the $\vc{a}$-bisect and $\vc{c}$-bisect are, respectively,
\begin{eqnarray*}
\vc{u}_a=\vc{a}-\vc{c}+w(\vc{c}-\vc{b}),  \\
\vc{u}_c=\vc{c}-\vc{b}+w(\vc{b}-\vc{a}),
\end{eqnarray*}
again with $w\in[0,1]$. Together the three collections of vectors give directions which cover a half circle. So for our given vector $\vc{u}$, exactly one of the equations in $w$ and $x$,
\begin{equation}\label{eqn:wxeqns}
\vc{u}x=\vc{u}_a,\ \ \vc{u}x=\vc{u}_b,\ \ \vc{u}x=\vc{u}_c
\end{equation} will solve for $w$ in $[0,1]$, unless $\vc{u}$ is in the direction of one of the triangle medians, in which case there will be two solutions, corresponding to $w=0$ and $w=1$. The sign of $x$ will indicate whether the orientation of $\vc{u}$ is consistent with our `right-to-left' definitions of $\vc{u}_a,\vc{u}_b$ and $\vc{u}_c$.

We return to these equations subsequently. But now we derive two vector equations from figure~\ref{fig:vectortriangle} which will allow us to solve for $t$, eliminating $s$. Without loss of generality we will again assume that $\vc{u}_b$ is the appropriate range of directions within which (up to orientation) lies vector $\vc{u}$. From figure~\ref{fig:vectortriangle} we see that, for the tilted line to form triangles with the edges $\vc{c}$ and $\vc{a}$ of $ABC$, the values of $s$ and $t$ must create resultant vectors in the direction of these edges. We express this as two orthogonality conditions:
\begin{eqnarray*}
 \left(-tm\vc{b}-s\vc{u}_b\right)\cdot \vo{c} &=& 0 \\
 \left((1-t)m\vc{b}+s'\vc{u}_b\right)\cdot \vo{a} &=& 0,
 \end{eqnarray*}
where $\vo{x}$, for a vector $\vc{x}$, denotes the perpendicular vector, whose dot product with $\vc{x}$ is zero.

Substituting for $\vc{u}_b$ and expanding out gives
\begin{eqnarray*}
 tm\vc{b}\cdot\vo{c}+s\vc{b}\cdot\vo{c} +s(w-1)\vc{a}\cdot\vo{c} &=& 0 \\
  (1-t)m\vc{b}\cdot\vo{a}+s'\vc{b}\cdot\vo{a}-s'w\vc{c}\cdot\vo{a} &=& 0,
 \end{eqnarray*}
 each of which may be factored using a single dot product by virtue of the fact that $\vc{a}+\vc{b}+\vc{c}=0$:
 \begin{eqnarray*}
 \left(tm+s  -s(w-1)\right)\vc{b}\cdot\vo{c} &=& 0 \\
  \left((1-t)m+s' +s'w\right)\vc{b}\cdot\vo{a} &=& 0.
 \end{eqnarray*}
 Since no triangle edge can be orthogonal to the perpendicular of another this gives
$$ tm+s(2-w) = 0, \hspace{.1in}
 (1-t)m+s'(1+w)= 0.$$
 We combine these two equations, using the substitution $s'=st/(1-t)$ from equation~\ref{eqn:sdash}, to arrive at $t^2/(1-t)^2=(2-w)/(1+w)$. Solving for $t$ gives
 $t=\left(1\pm\sqrt{(1+w)/(2-w)}\,\right)^{-1}$, from which we discard the negative square root since this gives a negative value of $t$ at $w=1$:
 \begin{equation}
 \label{eqn:tvalue}
 t=\left(1+\sqrt{\frac{1+w}{2-w}}\,\right)^{-1}.
 \end{equation}
The values of $t$, we observe, vary from $2-\sqrt{2}$ down to $\sqrt{2}-1$, as $w$ goes from~$0$ to~$1$. These extreme points mark precisely the intersections of the $\vc{a},\vc{b}$ and $\vc{c}$-bisects.

It remains to determine whether the value of $t$ identifies a point on the $\vc{b}$-bisect, as in figure~\ref{fig:vectortriangle}, or on one of the other two edge-bisects. To do this we must return to equations~\ref{eqn:wxeqns}. Solving these equations will identify an edge-bisect by giving a value of $w$ in $[0,1]$ which can be thereafter be substituted into equation~\ref{eqn:tvalue}. Now equations~\ref{eqn:wxeqns} are simultaneous equations in two variables, $x$ and $w$, which can be written down and solved explicitly by coordinatising $\vc{a}, \vc{b}$ and $\vc{c}$. The solutions may be expressed tidily by defining a function of three vector variables:
\begin{equation*}
P(\vc{x},\vc{y},\vc{z})=\frac{1}{\vc{x}\cdot \vo{y}}\matctwoone{\vc{z}\cdot \vo{y}}{\vc{z}\cdot \vo{x}},  \mbox{\ unless \ } \vc{x}\cdot \vo{y}= 0,  \mbox{\ in which case the function is undefined.}
\end{equation*}
The value of this function, if it is defined, is a column vector whose entries are $x$ and $w$. If it is not defined it is because $\vc{x}\cdot \vo{y}$ (which is the determinant of the matrix of our simultaneous equations) is zero, i.e. $\vc{x}$ is parallel to $\vc{y}$.

We proceed thus:
\begin{eqnarray}
\label{eqn:P}
&&\mbox{For\ }\vc{u}x=\vc{u}_a=\vc{a}-\vc{c}+w(\vc{c}-\vc{b}) \mbox{\ \ compute\ \ } P(\vc{u},\vc{c}-\vc{b},\vc{a}-\vc{c}),\nonumber\\
&&\mbox{For\ }\vc{u}x=\vc{u}_b=\vc{b}-\vc{a}+w(\vc{a}-\vc{c}) \mbox{\ \ compute\ \ } P(\vc{u},\vc{a}-\vc{c},\vc{b}-\vc{a}),\\
&&\mbox{For\ }\vc{u}x=\vc{u}_c=\vc{c}-\vc{b}+w(\vc{b}-\vc{a}) \mbox{\ \ compute\ \ } P(\vc{u},\vc{b}-\vc{a},\vc{c}-\vc{b}).\nonumber
\end{eqnarray}
Suppose that all three computations are defined, giving column vectors $(x_a,w_a)$, $(x_b,w_b)$ and $(x_c,w_c)$, say.
Then there will be a  unique value of $\{w_a,w_b,w_c\}$ lying in $[0,1]$. Let $t$ be the value of equation~\ref{eqn:tvalue} when $w$ takes this value. Then our required equation is
\begin{equation}\label{eqn:line}
    \vc{u}x+\begin{cases}
      \ \vcc{C}+(1-m)\vc{b}-tm\vc{a} & \text{if $w_a\in[0,1]$}  \\
      \ \vcc{A}+(1-m)\vc{c}-tm\vc{b}  & \text{if $w_b\in[0,1]$} \\
      \ \vcc{B}+(1-m)\vc{a}-tm\vc{c} & \text{if $w_c\in[0,1]$}
    \end{cases},\ \ x\in(-\infty,\infty).
\end{equation}
Suppose, on the other hand, that, say, $P(\vc{u},\vc{a}-\vc{c},\vc{b}-\vc{a})$ is undefined. Then $\vc{u}$ is parallel to $\vc{a}-\vc{c}$. But $\vc{a}-\vc{c}$ is the median from $B$ to the midpoint of $\vc{b}$. Then we may place $t$ on the $\vc{a}$-bisect, setting $w_a=0$, or on the $\vc{c}$-bisect setting $w_c=1$, and again specify the line equation from~\ref{eqn:line}.

We conclude with an example. Our triangle will be $A=(4, 2),\,B=(1, 9),\, C=(10, 1)$. The edge vectors are $\vc{a}=-B+C=(9, -8)$, $\vc{b}=-C+A=(-6, 1)$ and $\vc{c}=-A+B=(-3, 7)$, and the median (free) vectors are $\vc{b}-\vc{a}=(-15,9)$, $\vc{c}-\vc{b}=(3,6)$, and $\vc{a}-\vc{c}=(12,-15)$.

\begin{figure}[htb]
\centerline{\includegraphics[scale=0.6]{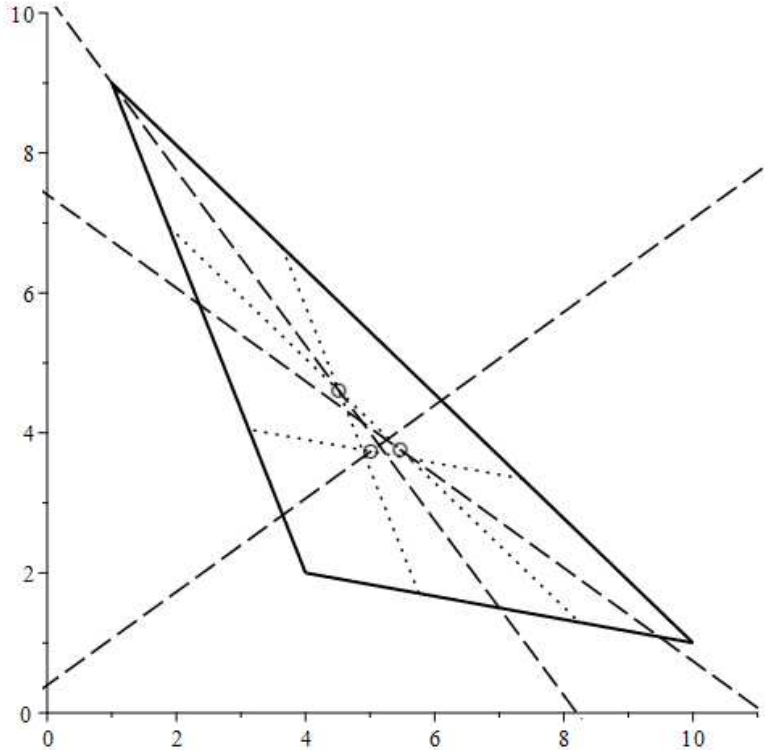}}

\caption{Triangle bisected by three lines (dashed) tilting on edge-bisects (dotted).}
\label{fig:example}
\end{figure}
In figure~\ref{fig:example} six bisecting lines are shown. The three dotted lines are the edge-bisects, the bisecting lines parallel to the triangle edges. The three dashed lines are calculated to achieve bisection for three different slope vectors $\vc{u}$, as follows:
\begin{description}
\item{$\vc{u}=(-3,2)$:} the three solution vectors $(x,w)$ computed at (\ref{eqn:P}) are, to two decimal places, $(-4.88,0.88)$, $(5.57,-0.14)$ and $(39.00,8.00)$, respectively. The first, $\vc{u}_a$ value of $w$ lies in $[0,1]$. Therefore the bisecting line passes through the point on the $\vc{a}$-bisect specified by $t$ (equation~\ref{eqn:tvalue}) at $w=0.88$, which is $t=0.44.$ This is plotted in figure~\ref{fig:example} as the dashed line meeting the vertical axis at approximately $7.4$.
\item{$\vc{u}=(3,2)$:} gives solution vectors  $(9.75,5.75)$, $(-1.70,0.83)$ and $(2.05,-0.21)$, respectively. This time $\vc{u}_b$ is selected. Again, to two decimal places, a $t$ value of $0.44$ is obtained and we may plot the dashed bisecting line meeting the vertical axis at approximately $0.4$.
\item{$\vc{u}=(-4,5)$:} gives solution vectors  $(-3,0)$, undefined, and $(3,1)$, respectively. Since the $\vc{u}_b$ computation is undefined, $\vc{u}=(-4,5)$ is parallel to the median from $B$ to the midpoint of $\vc{b}$. The other two computations give $w=0$ and $w=1$, as expected.
\end{description}

\end{document}